\documentclass[11pt]{amsart}
\usepackage{amsfonts,amsmath,amssymb}
\usepackage[all]{xy}
\usepackage{amsthm}

\usepackage[pagebackref]{hyperref}
\hypersetup{
	linktocpage,
	colorlinks = True,
	citecolor=blue,
	filecolor=blue,
	linkcolor=blue,
	urlcolor=black
}

\theoremstyle{definition}

\theoremstyle{plain}
\newtheorem{thm}{Theorem}

\begin{document}

\newcommand{\A}{{\mathcal A}}
\newcommand{\B}{{\mathbb B}}
\newcommand{\C}{{\mathbb C}}
\newcommand{\F}{{\mathbb F}}
\newcommand{\G}{{\mathbb G}}
\newcommand{\Gl}{{\rm Gl}}
\newcommand{\Sl}{{\rm Sl}}
\newcommand{\Z}{{\mathbb Z }}
\newcommand{\E}{{\mathcal E}}
\newcommand{\R}{{\mathbb R}}
\newcommand{\U}{{\rm U}}
\newcommand{\SU}{{\rm SU}}
\newcommand{\Spin}{{\rm Spin}}
\newcommand{\Sp}{{\rm Spec \;}}
\newcommand{\SO}{{\rm SO}}
\newcommand{\Oh}{{\rm O}}
\newcommand{\maps}{{\rm Maps}}
\newcommand{\three}{{\langle 3 \rangle}}
\newcommand{\sus}{{\Sigma^\infty}}
\newcommand{\ab}{{\rm ab}}
\newcommand{\Top}{{{\rm Top}_*}}
\newcommand{\ba}{{\bf a}}
\newcommand{\bc}{{\bf c}}
\newcommand{\be}{{\bf e}}
\newcommand{\ph}{{\hat{p}}}
\newcommand{\op}{{\rm op}}
\newcommand{\Aut}{{\rm Aut}}
\newcommand{\End}{{\rm End}}

\title{Toward a Galois theory of the integers over the sphere spectrum}

\author{Jonathan Beardsley and Jack Morava}

\address{The Johns Hopkins University and the University of Washington}

\email{jack@math.jhu.edu, jbeards1@uw.edu}

\thanks{The first author thanks Fred Cohen and Nitu Kitchloo, and the second author thanks Andrew Salch and Eric Peterson for helpful conversations. Both authors thank Oscar Randal-Williams for extensive tutorials about the matters in the appendix, and for allowing us to include an account of his work there. Both authors also thank the anonymous referee for many helpful comments.}

\begin{abstract}
	Recent work \cite{abg, abghr, beardsrelative,hess,rog} in higher algebra allows the reinterpretation of a classical description \cite{blumthh,cohenmaytaylor,mahinfinite,thomemss} of the Eilenberg-MacLane spectrum $H\Z$ as a Thom spectrum, in terms of a kind of derived Galois theory. This essentially expository talk summarizes  some of this work, and suggests an interpretation in terms of configuration spaces and monoidal functors on them, with some analogies to a topological field theory. 
\end{abstract}

\maketitle

\section{Introduction}

\subsection{} Questions about an absolute base $\F_1 \to \Z$ for arithmetic (e.g.~\cite{connesconsanif1,kapranovsmirnov}) evoke echoes in algebraic topology (e.g.~\cite{toenvaquie,toenvezzosiringspectra}), where Waldhausen's brave new rings program interprets the unit $1 \in \tilde{H}^0(S^0,\Z)$ as a kind of ring homomorphism from the sphere spectrum $S^0$ to the Eilenberg-MacLane spectrum $H\Z$. This talk, aimed at interested non-experts, tries to summarize current thinking (and speculation) about this and related questions in homotopy theory, framed in terms of recent work of Rognes, Hess, and others on an emerging version of Galois theory in higher algebra. A now classic construction by Mahowald from the late 1970s (then at $p = 2$, though quickly generalized \cite[Corollary 3.5]{mahinfinite} \cite{cohenmaytaylor}) interprets $H\Z$ as a Thom spectrum (or cobordism theory); this talk describes that result in this developing language. 

{\bf Some conventions:} We write $\wedge$ for the (symmetric monoidal) smash product of pointed spaces, and $\wedge_R$ (e.g.~$\wedge_{S^0}$) for the smash product of $R$-module spectra. If $G$ is a group, the geometric realization or nerve $|[*/G]| := BG$ of the associated singleton category provides a standard model for its classifying space. $\sus X_+$  (or $S^0[X]$) will denote the suspension spectrum defined by an unpointed space. $X_\infty$ will denote one-point compactification (e.g.~$\C_\infty \cong S^2$). For a space $X$ we use $X\langle n\rangle$ to indicate the $n$th level of the Whitehead tower, i.e.~the space which has trivial homotopy in degrees less than \emph{or equal to} $n$, and is equivalent to $X$ elsewhere. There seems to be some disagreement in the literature about the precise meaning of this notation and how it relates to the terms ``$n$-connected'' and ``$n$-connective,'' but we will always mean the space described above.

\subsection{} The homotopy theorist's category of spectra (cf.~\cite{ekmm,symmspectra}), or modules over the sphere spectrum, is the distinguished ancestor of the modern theory of derived categories and their relatives (such as Voevodsky's motives), and for our principally expository purposes we will assume some general familiarity with its constructions. However, a significant difference between classical and homotopy-theoretic ring objects is the subtlety of notions of associativity and commutativity in the latter context. There, a multiplicative structure on an object $A$ can usefully be presented by a compatible collection
\[
\E_n(k) \wedge_{\Sigma_k} A^{\wedge k} \to A
\]
of parametrized multiplication morphisms, making $A$ an algebra over an {\bf operad} $\E_n := \{\E_n(k), \; k \geq 0\}$  of objects with an action of the system $\{\Sigma_k\}$ of symmetric groups. In our situation $\E_n(k)$ will be, very roughly, a space of ordered $k$-tuples in $\R^n, \; 1 \leq n \leq \infty$ \cite{mayquinnray}. An $\E_1$ (or $A_\infty$) algebra structure thus defines a coherently homotopy-associative product, and an $\E_\infty$ algebra is similarly coherently homotopy-commutative. In particular, the unit map
\[
S^0 \to H\Z
\]
is a morphism of $\E_\infty$ algebras (in the category of spectra, i.e.~$S^0$-modules), but the description of $H\Z$ as a Thom spectrum presents it only as an $\E_2$-algebra. 

When $n>2$ the spaces $\E_n(k)$ are simply-connected, but $\E_2(k)$ is a classifying space for the $k$th pure braid group: $\E_2$-algebras thus have a {\bf braided} monoidal structure. This suggests possible connections with the study of automorphisms of the braid groups (i.e.~understanding the work of Grothendieck and Drinfel'd on the tower of mapping class groups  \cite{bhrGT,horelprof}); we won't discuss this further here, but it is one of the motivations for our interest in this subject \cite{kitchloomorava,moravamotives}.

\subsection{} The organization of the paper is as follows: in Section \ref{section:1hges} we briefly recall the notion of a Hopf-Galois extension, for both fields and ring spectra. We also state our main theorem, which is that the unit map $S^0\to H\Z$, when considered as a map of $\E_2$-algebras, is such an extension. This can be thought of as a description of the descent data for passing from algebra, over $\Z$, to homotopy theory, over $S^0$. 

Our main theorem can be considered a restatement of the material we review in Section \ref{section:2thomspectra}. Therein, we recall the important work of Ando, Blumberg, Gepner, Hopkins and Rezk in describing Thom spectra as quotients of actions by grouplike iterated loop spaces. We also recall the descriptions of $H\F_p$, $H\Z_\ph$ and $H\Z$ as Thom spectra, described by Blumberg as well as Antolin-Camarena and Barthel. In particular, $H\Z$ is described as the Thom spectrum of a 2-fold loop map $\Omega^2(S^3\langle 3\rangle)\to BO$. 

Finally in Section \ref{section:3braidmodel}, which we include as an appendix, we give a concrete categorical model for the space $\Omega^2(S^3\langle 3\rangle)$ in terms of braids of writhe, or total twist, equal to zero.

\section{Hopf-Galois theory in Spectra}\label{section:1hges}

\subsection{} The normal basis theorem of elementary Galois theory can be phrased as the existence, for a finite normal field extension $F \subset E$ with Galois group $G$, of an isomorphism
\[
E \otimes_F E \cong F_G \otimes_F E
\]
of $F$-vector spaces \cite{kreitak}: where $F_G$ is the Hopf-algebra of functions from $G$ to $F$, dual to the group ring $F[G]$. This realization allows one to generalize the notion of a $G$-Galois extension to extensions of fields (or rings, as in \cite{kreimergalois}) whose ``Galois groups'' are in fact Hopf-algebras or even just bialgebras \cite{lebruynF1,montsurv}. We will use the terminology \textbf{Hopf-Galois object} for the generalized Galois groups of such extensions. 

About ten years ago Rognes \cite{rog} took this as the basis for a Galois theory of $\E_\infty$ ring morphisms in the category of spectra, arguing in particular that the Thom isomorphism
\[
M\U \wedge_{S^0} M\U \simeq (\sus B\U_+) \wedge_{S^0} M\U
\]
for the complex cobordism spectrum $M\U$ (with the suspension spectrum $\sus B\U_+$ regarded as a strictly cocommutative bialgebra with diagonal induced by $B\U_+ \to B\U_+ \wedge B\U_+$, and multiplication defined by the ($\E_\infty$) Whitney sum map $\oplus : B\U_+ \wedge B\U_+ \to B\U_+$) can be interpreted as making the unit map
\[
S^0 \to M\U
\]
an $\E_\infty$ extension of ring spectra, with $\sus B\U_+$ as its Hopf-Galois object \cite[\S 4.1.3]{rog}. [Note that $B\U_+$ is not strictly dualizeable; if it were, $\sus B\U_+$ might be interpreted as a Hopf algebra of functions on some Spanier-Whitehead dual ``spectral groupscheme.'']

 Applications and examples of these ideas followed quickly: Hess \cite{hess} used them as the basis for a derived Tannakian theory of homotopical descent and codescent; others \cite{froth} considered weakening the $\E_\infty$ conditions, and developments in higher category theory (discussed below) led to a theory of intermediate extensions interpreting, for example, the forgetful morphism from complex to oriented cobordism, $MU\to MSO$, as Galois, with $\sus \Spin_+$ as Hopf-Galois object \cite[\S 2]{beardsrelative}. Other connections with the theory of Brauer and Picard groups will be relevant below.

\subsection{} It is classical \cite{atiyahhirz,milnorsteenrodalg} that the Hopf algebra dual  
\[
\pi_* (H\F_2 \wedge_{S^0} H\F_2) \cong \A_* = \F_2[\xi_i \: | \: i \geq 1]
\]
to the Steenrod algebra is polynomial, on generators $\xi_i$ of degree ${2^i - 1}$. Its spectrum, in the sense of commutative algebra, is thus the dual $({\rm Prim} \; \G_a)^\vee$ of the space of primitive elements of the additive groupscheme $\G_a := \Sp \F_2[T]$, and represents the functor of endomorphisms of $\G_a$.  It is similarly classical that the Pontrjagin (Hopf) algebra 
\[
H_*(\Omega^2 S^3,\F_2) \cong \F_2[x_i \:|\: i \geq 1] 
\]
($|x_i| = 2^i - 1$) is primitively generated, and is in fact free on one generator over the Kudo-Araki operations: $x_i = \bar{Q}^i x_1$ (where $\bar{Q}^i$ is the iterated Kudo-Araki operation $Q^{2^{i-1}}\circ Q^{2^{i-2}}\circ\cdots\circ Q^2$); and $x_1$ is the Hurewicz image of the adjoint
\[
S^1 \to \Omega^2 S^3 
\]
to the identity map of the three-sphere \cite[Theorem 7.1]{kudoaraki}. Moreover, the cohomology $H^*(\Omega^2 S^3,\F_2)$ is free of rank one over the Steenrod algebra. 

Mahowald's construction interprets this as an equivalence of $\E_2$-$H\F_2$-algebras
\[
H\F_2 \wedge_{S^0} H\F_2 \simeq (\sus \Omega^2 S^3_+) \wedge_{S^0} H\F_2 
\]
coming from a presentation of $H\F_2$ as a Thom spectrum defined by the map
\[
\Omega^2 S^3 \to B\Oh.
\]
This map is obtained from
\[
[\eta : S^1 \to B\Oh] \in KO_1 \cong \{\pm 1 \}
\]
using Bott's infinite-loopspace structure on $B\Oh$: a map from $X$ to an $n$-fold loopspace $Y$ extends, via the diagram
\[
\xymatrix{
X \ar[d] \ar[r]^h & Y \simeq \Omega^n Z \\
\Omega^n S^n X \ar@{.>}[ur]^{\alpha(h)} \ar[r] & \Omega^n S^n \Omega^n Z 
\ar[u] \;,}
\]
to an $n$-fold loop map $\Omega^n S^n X \to Y$. In terms of $\E_2$-algebras, this construction comes from an equivalence 
\[
\maps_\Top(S^1,B\Oh) \simeq \maps_{\E_2 - {\rm alg}}(\Omega^2 S^3,B\Oh)
\;. \]
Comparing the Dyer-Lashof and Steenrod algebra structures \cite[\S 9]{blumthh} implies that 
\[
S^0 \to H\F_2
\]
is an extension of $\E_2$-algebras with $\sus \Omega^2 S^3_+$ as Hopf-Galois object.

\subsection{} The $\F_1$-approach to the Riemann hypothesis proposes $\Sp \Z  \times_{\F_1} \Z$, or a suitable compactification, as an analog of Weil's space $C \times_k C$ for a curve $C$ over a finite field $k$: with an associated algebra of correspondences built from its divisors \cite[\S 2.1]{ccm}. In an attempt to understand $\mathbb{Z}$ as a curve over $S^0$, we could try to understand the derived algebro-geometric properties of the tensor product $H\Z\wedge_{S^0}H\Z$, or of the functor of points defined by its homotopy groups.

In classical algebraic topology, the dual Steenrod algebra $\mathcal{A}_\ast$ corepresents the groupscheme $\Aut_{\hat{\G}_a}$ of automorphisms of the additive formal group scheme $\hat{\G}_a$ (cf. \cite[\S 1]{inouesteenrod} or \cite{milnorsteenrodalg}). This functor is given by
\[
\F_2 - {\rm alg} \ni A \mapsto \left\{a(T) = T + \sum_{i \geq 1}a_i T^{2^i} \in 
A[[T]]\right\} := \Aut_{\hat{\G}_a}(A)
\]
and its elements can explicitly seen to be automorphisms,
\[
a(T_0 + T_1) = a(T_0) + a(T_1),
\]
of the additive formal group $\mathrm{Spf}(A[[T]])$ over $A$, where Spf is the formal spectrum functor. The groupscheme structure here comes from composition of power series. If we think of the $\F_2$-homology $H^\ast(X,\mathbb{F}_2):=\pi_\ast(\Sigma^\infty X_+\wedge_{S^0}H\F_2)$ as a kind of fiber functor 
\[
({\rm finite \; spectra}) \ni X \mapsto H^*(X,\F_2) \in (\F_2 - {\rm Vect}) \;.
\]
then we can think of the groupscheme $\Aut_{\hat{\G}_a}$ as (the dual of) the automorphisms of this functor, in the Tannakian formalism of \cite[Chapter II]{milnedelignetannaka}.

More generally, by identifying $\A_*$ with the Hopf-algebra endomorphisms of the (usual, not formal) additive group $\G_a$, the Steenrod algebra can be interpreted as corepresenting
the $\F_2$-algebra-valued functor
\[
A \mapsto A \langle F \rangle := \left\{ \sum_{n \geq 0} a_n F^n \:|\: i \gg 0 
\Rightarrow  a_i = 0\right\}
\]
(with semilinear multiplication ($Fa = a^2F, \; a \in A$)), where $F$ is the Frobenius endomorphism, see also \cite{smithalgebraicsteenrod}. This suggests regarding $\End_{\G_a}$ as an analog, in the category of ring-schemes, of an algebra of correspondences for the integers over $\F_1$. This algebra is in some sense of infinite rank, but the $\F_1$-program envisages $\overline{{\rm Spec} \; \Z}$ as a curve of infinite genus anyway. 

\subsection{}
Following \cite{hess}, we may also think of $H\F_2\wedge_{S^0}H\F_2$ as the \textbf{descent coring} of the $\E_\infty$-morphism $S^0\to H\F_2$ controlling descent from $\F_2$-vector spaces to $S^0$-modules. On the other hand, Mahowald's construction of $H\F_2$ as a Thom spectrum indicates that if we only concern ourselves with $S^0\to H\F_2$ as an $\E_2$-map, we can replace this descent coring with $\Sigma^\infty\Omega^2 S^3_+ \wedge_{S^0}H\F_2$. Moreover, recall from e.g.~\cite{montsurv} that the Hopf-Galois extension associated to a $G$-Galois extension will have Hopf-algebra the (dual of the) ring of functions on $G$. This suggests that we may think of the primitively generated Hopf algebra $H_*(\Omega^2S^3,\F_2)$ as the symmetric algebra of functions on the affine space underlying $\End_{\G_a}$. Its generators $x_n$ then become conceptually analogous to the divisors generating Weil's algebra of correspondences. If, in light of \cite{balmspec, balmerkrausestevenson}, we consider the Morava K-theories $K(n)$ to be the prime ideals of the sphere spectrum, then the generators $x_n\in H_*(\Omega^2S^3,\F_2)$ correspond to these points thought of as divisors in the fiber product $Spec(H\F_2)\times_{Spec(S^0)}Spec(H\F_2)$. [At the moment, a structure space or spectrum for an $\E_2$-ring spectrum is yet to be defined, but the convergence of the classical Adams spectral sequence reflects the fact that this fiber product ``knows'' about the primes of $S^0$.] 

One would like there to be a similar Hopf-Galois object for the unit maps $S^0\to H\F_p$ when $p\neq 2$ as well. Unfortunately, the isomorphism $\pi_1B\Oh \cong \Z/2$ precludes the possibility of a straightforward generalization of this construction to odd primes (in particular, for $p>2$, it is not the case that $p-1$ is invertible in $\pi_0(S^0)\cong\mathbb{Z}$). Any Thom spectrum arising from a map $X\to B\Gl_1(S^0)$ must have $\pi_0$ isomorphic to either $\Z$ or $\Z/2$. However, following ideas of Hopkins (as described in \cite[Lemma 3.3]{thomemss}, \cite[\S 9.2]{blumthh} and \cite[Theorem 5.1]{acb}), one obtains the same results at odd primes by working over the $p$-localization of $B\Gl_1(S^0)$. This kind of construction has also recently been extended to include $H\Z_{p^k}$ in \cite{kitchloothh}.  

We are, however, ultimately interested in an integral extension of $S^0$ which relates the global structure of $\Z$ to that of $S^0$ via a kind of descent from algebra to topology. This goal is reflected in the following:

\begin{thm}
The unit morphism 
\[
S^0 \to H\Z
\]
is an $\E_2$ - Galois extension with $\sus \Omega^2(S^3\three)_+$ as Hopf-Galois object.
\end{thm}

Here $S^3\three$ is the three-connected cover of 
$S^3$  (which is of course homeomorphic to $\SU(2)$; also note that $\SU(2)\three$ is known in some physics circles as $\mathrm{String}(1)$). In proving this statement, one first proves its $p$-adic analogues; that is, one proves equivalences $$H\Z_\ph\wedge_{S^0_\ph} H\Z_\ph\simeq H\Z_\ph\wedge_{S^0_\ph} \sus \Omega^2(S^3\three_\ph)_+.$$ Here the subscript $\ph$ signifies $p$-adic completion, so $\Z_\ph$ is the 
ring of $p$-adic integers. This in turn follows from the construction of $H\Z_\ph$ as a Thom spectrum over a 2-fold loop space \cite[Theorem 5.7]{acb} \cite[\S 9.3]{blumthh}, as the above equivalence is precisely the Thom isomorphism. 

Local-to-global methods \cite[Lemma 9.3]{blumthh} then lead to a description 
of $H\Z$ as a Thom spectrum over the double loop space $\Omega^2(S^3\three)$, and thus an $\E_2$-extension of $S^0$ with Hopf-Galois object 
$\sus \Omega^2 (S^3\three)_+$.

\section{Thom spectra from spherical fibrations}\label{section:2thomspectra}

\subsection{} The Thom space of a vector bundle $E \to X$ is the cofiber of the projection $S(E) \to X$ of its unit sphere bundle; this suggests the extension of the theory of Thom spectra to more general spherical fibrations, classified not by maps to $B\Oh$ but to a classifying space $B\Gl_1(S^0)$ for the monoid of stabilized self-equivalences of the sphere. More generally, Ando, Blumberg, Gepner, Hopkins, and Rezk \cite{abg,abghr} use $\infty$-categorical methods to associate to an $\E_1$-algebra spectrum $R$, a classifying space $B\Gl_1(R)$ for ``bundles'' of free rank one$R$-modules (which they call $R$-lines); if $R$ is an $\E_\infty$ - algebra, this classifying space is an infinite loop space. To be precise, $B\Gl_1(R)$ is the Kan complex or $\infty$-groupoid of $R$-lines and equivalences between them, hence the base-point component of the so-called \textbf{Picard category} of $\mathrm{LMod}_R$. So defined, there is a natural inclusion of $\infty$-categories $B\Gl_1(R)\hookrightarrow \mathrm{LMod}_R$. By regarding a space X as an $\infty$-category in \cite[\S2.5]{abghr}, these authors  associate to a 
map $f : X \to B\Gl_1(R)$, a Thom $R$-module spectrum 
\[
\xymatrix{
Mf := {\rm colim} ( X \ar[r]^-f & B\Gl_1(R) \ar@{^{(}->}[r] & {\rm LMod}_R
 ) .} 
\]
Moreover, they show that if $G$ is a group-like monoidal $\infty$-groupoid, then a map $BG \to B\Gl_1(R)$ defines an $R$-linear action of $G$ on $R$, yielding an equivalence (cf.~\cite[Theorem 1.17]{abghr}):
\[
Mf \simeq R//G 
\]
of their Thom spectrum with the homotopy quotient (n.b.~not the homotopy fixed-points: this is not our Grandmother Emma's Galois theory!) of $R$ by $G$. For example, the classical $J$-homomorphism defines a map
\[
B\Oh \to B\Gl_1(S^0)
\]
and thus an identification $M\Oh \simeq S^0//\Oh$.\bigskip

Hopkins' reinterpretation of Mahowald's theorem is based on an identification
\[
\pi_1 B\Gl_1(S^0_\ph) \cong \Z^\times_\ph 
\]
of the fundamental group of the Picard category of $S^0_\ph$-lines with the multiplicative group of invertible $p$-adic integers. A map $h : S^1 \to B\Gl_1(S^0_\ph)$ extends, as in \S 2.2, to a two-fold loop map $\alpha(h) : \Omega^2 S^3 \to B\Gl_1(S^0_\ph)$, defining a homomorphism $\pi_1\alpha(h) : \Z \to \Z^\times_\ph$ on fundamental groups. If, for example, $h = 1 - pk$ with $(p,k) = 1$ then $\pi_1\alpha(h)$ maps $1 \in \Z$ to a topological generator of $(1 + p\Z_\ph)^\times$, defining an isomorphism of $\pi_1(\Omega^2 S^3)$ with a dense subgroup of $\pi_1B\Gl_1(S^0_\ph)$. The Thom spectrum of this morphism is then $H\F_p$ \cite{acb,blumthh,thomemss}. Lifting the map $\alpha(h)$ to universal covers defines a double loop map
\[
\Omega^2(S^3\three) \simeq B(\Omega^3_0S^3) \to  B\Sl_1(S^0_\ph),
\]
whose Thom spectrum is $H\Z_\ph$, with $\Omega^3_0S^3$ the base-point component of $\Omega^3S^2$ \cite[Theorem 5.7]{acb}. From \cite{blumthh} we then have that the $p$-adic maps may be glued together to produce $H\Z$ as a Thom spectrum over $\Omega^2(S^3\three)$.  This defines an equivalence of $\E_2$-algebras:
\[
H\Z \simeq S^0//\Omega^3_0S^3 \;.
\]

It will be relevant below that the image 
\[
\pi_{2(p-1)}\alpha(h) : \pi_{2(p-1)}\Omega^2(S^3\three) = \pi_{2p}S^3 \to \pi_{2(p-1)}B \Gl_1(S^0) = \pi^S_{2p-3}
\]
contains the class $\alpha_1$ (generating the image of the $J$-homomorphism in that degree \cite[Theorem 3]{cohenmaytaylor}). \bigskip

\subsection{}\label{braidrep} At this point it is interesting to reconsider the $p=2$ case, and recall that 
\[
\alpha(h) : \Omega^2 S^3 \to B\Oh (\simeq \Omega^2 {\rm Sp})
\]
can be constructed by group-completing the following sequence of monoids:
\[
\coprod_{n \geq 0} B\B_n  \to \coprod_{n \geq 0} B\Sigma_n \to \coprod_{n \geq 0}
B\Oh(n)
\] 
Here $\B_n$ is the $n$-strand braid group (discussed in more detail in the appendix) and $\Sigma_n$ is the symmetric group on $n$ elements. The map $\B_n\to \Sigma_n$ is the ``underlying permutation'' map that takes a braid to the permutation obtained by forgetting all crossings and only remembering the endpoints of the braids. The map $B\Sigma_n\to B\Oh(n)$ is the defined by the regular representation which permutes coordinates. Here, following \cite{segalconfig}, $B\B_n$ can be interpreted as a configuration space of $n$ points in $\C$, so the composite
\[
\coprod_{n \geq 0} B\B_n \to \Omega^2 S^2 \to \Z \times B\Oh
\]
can be regarded as a kind of topological field theory $F \mapsto \R^F$ which takes a configuration $F$, with $\#F=n$, to $\R^n$, and a braid from  $F$ to $F'$ to the transformation $\R^n\to\R^n$ which permutes coordinates in accordance with the forgetful map $\B_n\to\Sigma_n$. More precisely, the domain category has finite subsets (codimension two submanifolds) of $\C$ as its objects\begin{footnote}{Perhaps with an underlying $\mathbb{F}_1$-module?}\end{footnote}, and braids, regarded as isotopy classes of one-dimensional cobordisms (codimension two submanifolds embedded in $[0,1] \times \C$) as morphisms; while the range is the category $\mathrm{Vect}_\mathbb{R}$ with monoidal structure given by $\oplus$. The monoidal structure in the domain, which is given by juxtaposition of braids, can be defined geometrically on configurations in $\mathbb{C}\subset \mathbb{C}_\infty$ by pulling back along a suitable pinch map $\C_\infty \to \C_\infty \vee \C_\infty$, e.g.~constructed by collapsing $\R_\infty \subset \C_\infty$ to a point. This suggests asking if something analogous might underlie the Thom space interpretation of $H\Z$, defined by a functor from some such configuration space to a Picard category of $S^0$-lines.

In fact classic work on `finite models' for the $J$-homomorphism, going back to \cite[Chapter 2]{madmil} \cite{tsuchiya} may be relevant to this question. Sullivan showed that the $J$-homomorphism splits as a map of spaces, and later work (e.g.~\cite{friedlanderadams}) extends this to a splitting (at least at odd primes): 
\[
\Gl_1(S^0) \simeq J \times {\rm Coker \;} J
\]
as infinite loopspaces. In more modern language we might think of this splitting as defined by Rezk's logarithm \cite{rezklog}
\[
R\log_K : {\rm gl}_1(S^0) \to L_KS^0
\]
(cf.~e.g.~\cite[\S 8]{rav3}) in one direction, where $\mathrm{gl}_1(S^0)$ is the spectrum associated to the infinite loopspace $\Gl_1(S^0)$, and in the other by Tornehave's equivalence
\[
J^\oplus_p \simeq J^\otimes_p
\]
of infinite loopspaces (cf.~\cite[VIII Corollary 4.2, Remark 4.6]{mayquinnray} and \cite{mayJsplit}), based on Quillen's models for the image of $J$ in terms of $\F_l$-vector spaces ($l \equiv 1$ mod $p^2$), viewed alternatively as finite sets under Cartesian product \cite{mstgtop,snaithloopmaps}. The issue is that the monoidal structure on $\Omega^3_0S^3$ comes from loop addition, while that of $\Gl_1(S^0)$ comes from smash product; but (at odd primes) the morphism 
\[
\Omega^3_0 S^3_\ph \to \Sl_1(S^0_\ph)
\]
factors through the image of $J$, where the distinction between $\oplus$ and $\otimes$ simplifies. Note that by $\Sl_1(R)$, for a ring spectrum $R$, we mean the connected component of the identity in $\Gl_1(R)$. 

This can perhaps be summarized in a diagram of the form
\[
\xymatrix{
S^1 \ar[d] \ar[r] & B\Gl_1(S^0_\ph) \ar[r]^-R & \ar[l]^-T BJ_p \\
\Omega^2 S^3 \ar[ur]^{\alpha(h)} & B\Sl_1(S^0_\ph) \ar[u] \ar[r]^-R & 
\ar[l]^-T B\tilde{J}_p \ar[u]\\
\Omega^2(S^3 \three) \ar[u] \ar[ur]^{\tilde{\alpha}(h)} \;,}
\]
interpreted as defining a kind of monoidal functor from a suitable category of ``configurations'' to $L_KS^0$-lines where $R$ and $T$ are the maps of Rezk and Tornehave, respectively. At odd primes $p$ we may think of  replacing $L_KS^0$ by ($\F_l$-Mod,$\otimes$). The appendix below describes a tentative model for such a category of configurations.

\section{Appendix: A categorification of $\Omega^2 (S^3 \three)$}\label{section:3braidmodel}

\subsection{} We show in what follows that the fibration $\Omega^2S^3\langle 3\rangle\to \Omega^2S^3\to S^1$, as a fibration of 2-fold loop spaces, can be constructed from a fibration whose components are group completed classifying spaces of categories. The space corresponding to $\Omega^2S^3$ is the free braided monoidal category on an object $\B_{\bullet}$, $S^1$ arises from the abelianization of $\B_\bullet$, and the fiber arises from a category of braids of zero writhe.  

 Let 
\[
\B_\bullet := \coprod_{n \geq 0} [\{n\}/\B_n]
\]
be the free braided monoidal groupoid on one generator as in \cite{joyalstreet}, i.e.~with the set $\Z_{\geq 0}$ of non-negative integers as objects, and morphisms:

\[
\maps_{\B_\bullet}(m,n)=\left\{ \begin{matrix}
	\{1\} & n=m\leq 1\\
	\B_n & n=m>1\\
	\emptyset & n\neq m\\
	\end{matrix}\right.
\]

where $\B_n$ is Artin's braid group   
\begin{align*}
\B_n =\langle \sigma_i, 1 \leq i \leq n-1 \:|\: \sigma_i \sigma_{i+1} \sigma_i = 
\sigma_{i+1} \sigma_i \sigma_{i+1}, \; |i-j| > 1 \Rightarrow 
[\sigma_i,\sigma_j] = 1 \rangle.
\end{align*}
Its monoidal structure is defined on objects by the functor $(n,m)\mapsto n+m$ and on morphisms by the juxtaposition product
\[
\B_n \times \B_m \to \B_{n+m}
\]
that constructs an $n+m$ stranded braid by placing an $n$-stranded braid next to an $m$-stranded one. Let $c_{n,m}$ be the braid that passes the trivial $n$-strand braid on the left over the trivial $m$-strand braid on the right. Then $c_{n,m}:n+m\to m+n$ defines a braiding on $\B_\bullet$ \cite[Example 2.1]{joyalstreet}.  

\subsection{}\label{abelianizationmap} If $n > 1$, the abelianization map
\[
w_n : \B_n \ni \sigma_i \mapsto 1 \in \Z \cong \B^\ab_n
\]
sends a braid to its total twist or {\bf writhe}, i.e.~the number of overcrossings less the number of undercrossings; if $n = 0$ or 1 we take $w_n$ to be $0 \to \Z$. Being a groupoid, $\B_\bullet$ has an abelianization (constructed by abelianizing morphism objects): 
\[
\B_\bullet \to \B^\ab_\bullet := [\Z_{n \geq 0}/\Z]
\]
(defined by the trivial action of the integers on $\Z_{n \geq 0}$) which is again monoidal, with a nontrivial braiding. This functor is the identity on objects of $\B_\bullet$ and on mapping objects takes each braid to its writhe. We are very grateful to Oscar Randal-Williams for explaining the properties of this object to us, and for his permission to include parts of his description here. Of course any mistakes or falsities are wholly a result of our own misunderstanding! 

The monoidal structure on $\B^\ab_\bullet$ is given by $(p,q)\mapsto p+q$ on objects and by addition on maps:
\[
\maps_{\B^\ab_\bullet}(p,q)\times \maps_{\B^\ab_\bullet}(r,s)=\Z\times\Z\overset{+}{\to} \Z = \maps_{\B^\ab_\bullet}(p+r,q+s).
\]
It is clearly associative and unital, with associator given by associativity of addition and unitality given by the fact that $a+0=a$. There is in fact a symmetric monoidal structure on this category, but the monoidal structure induced by the abelianization functor is only braided. Recall that the braiding in $\B_\bullet$ was just the invertible action of $\B_{n+m}=\B_{m+n}$ on itself by left multiplication with the Joyal-Street element $c_{n,m}$. The element $c_{n,m}$ abelianizes to 
\[
nm \in\Z = \maps_{\B^\ab_\bullet}(n+m,m+n),
\]
(since it has $n$ strands crossing over $m$ strands). In other words, the braiding of $\B_\bullet^\ab$ is the isomorphism $n+m\overset{nm}\longrightarrow n+m=m+n$ whose inverse is $-nm$.  The groupoid abelianization functor preserves monoidal structure, so it is immediate that this defines a braiding on $\B_\bullet^\ab$. For the sake of clarity, we verify the relations defining a braided monoidal structure on $\B_\bullet^\ab$. 

Let $p$, $q$ and $r$ be objects of $\B_\bullet^\ab$ (hence non-negative integers). Then, again following \cite{joyalstreet}, we must have commutativity of the following diagrams (where $0$, the identity map, is the associator):

\[
\xymatrix{
	(p+q)+r \ar[d]^0 \ar[r]^{pq} & (q+p)+r \ar[d]^{(q+p)r} \\
	p+(q+r) \ar[d]^{qr}  & r+(q+p) \ar[d]^{0} \\
	p+(r+q) \ar[r]^{p(r+q)} & (r+q)+p,}\hspace{.6in}
\xymatrix{
	p+(q+r) \ar[d]^0 \ar[r]^{qr} & p+(r+q) \ar[d]^{0} \\
	(p+q)+r \ar[d]^{(p+q)r}  & (p+r)+q \ar[d]^{pr} \\
	r+(p+q) \ar[r]^{0} & (r+p)+q.}
\]

Taking geometric realization of these categories defines a map
\[
|\B_\bullet| \to |\B^\ab_\bullet|
\]
of topological monoids: in fact of algebras over the $\E_2$ operad (one can show by an elementary computation with Browder brackets that these cannot extend to $\E_3$-algebra structures). We would like to understand the homotopy type of this map and its fiber, after group completion. To that end, we first recall the homotopy type of the source. 

Classical work \cite{boardvogt,segal} of Boardman, Vogt, and Segal describes a homotopy equivalence of $|\B_\bullet|$ with the space of finite subsets of $\C$
\[
\{ F \subset \C \:|\: \# F < \infty\},
\]
of finite subsets of the plane, topologized as a coproduct of configuration spaces $Conf_{n}(\C)$. Segal further points out that to each configuration, by thinking of the points as ``electrical particles with charge $+1$,'' one may associate an electrical field $E_c:\C\setminus F\to \C$. This extends to a map $E_c:\C_\infty\to \C_\infty$ by $E_c(p)=\infty$ whenever $p\in F\cup\{\infty\}$. Concretely, we may either think of this construction as a sum of rational functions with poles at $F$, or by the Weierstrass Factorization Theorem, as a polynomial with zeros at $F$: 
\[
|\B_\bullet| \ni F \mapsto p_F(z) = \prod_{\rho \in F}(z - \rho) \in \C[z] \;.
\]
In either case we obtain a function of topological spaces that Segal calls the ``completion map''
\[
|\B_\bullet|\to\maps_\Top(\C_\infty,\C_\infty)\simeq \Omega^2 S^2.
\]
Segal further shows (cf.~also \cite{maygeom}) that this map is an equivalence after delooping once. Thus, by thinking of the functor $\Omega B$ as \textbf{group completion}, we have an equivalence: 
\[
\Omega B|\B_\bullet|\overset{\sim}\to \Omega B \Omega^2S^2 \simeq \Omega^2S^2.
\]

Now note that we have $|\B_{\bullet}^\ab|=|\coprod_{n\geq 0} [\{n\}/\Z]|\simeq \mathbb{N}\times S^1.$ Hence group completion induces an equivalence $\Omega B|\B_{\bullet}^\ab|\simeq \Z\times S^1$. Recalling that $\Omega^2 S^2\simeq \Z\times\Omega^2S^3$,  it follows that the map 
\[\Omega B|\B_{\bullet}|\to \Omega B|\B_{\bullet}^\ab|\]
 is the 1-type truncation 
 \[\Omega^2S^2\simeq \Z\times\Omega^2S^3\to \Z\times S^1\simeq \Omega^2S^2[0,1]\simeq \Omega^2(S^2[0,3]).\] 
 Here $X[0,n]$ denotes the space obtained by killing all homotopy groups of $X$ above $\pi_n$. 
 
 Because this map is obtained as the group completion of the geometric realization of a braided monoidal functor, we know that it is a map of $\E_2$-algebras, but we can give a more precise construction. In particular, note that the 3-type truncation $S^2[0,3]$ is defined by a Postnikov fibration 
 \[
 S^2[0,3] \to K(\Z,2) \to K(\Z,4),
 \]
 where the second map is $x^2\in H^\ast(K(\Z,2);\Z)\cong \Z[x]$. By taking the two-fold delooping of the truncation map $S^2\to S^2[0,3]$ we obtain our 2-fold loop map 
 \[
 \Omega^2S^2\to \Omega^2(S^2[0,3]).
 \] 
Moreover, noticing that the attaching map $x^2:K(\Z,2)\to K(\Z,4)$ must be trivial after looping once (because $K(\Z,1)\simeq S^1$ has no cohomology in degree 3), we obtain an $\E_2$ splitting $\Omega^2(S^2[0,3])\simeq \Z\times S^1$. Finally, the fiber of the identity component of the map 
\[
\Omega B |\B_\bullet| \simeq \Omega^2 S^2 \simeq \Omega^2 S^3 \times \Z 
\to \Omega B |\B^\ab_\bullet| \simeq \Omega^2 (S^2[0,3])
\]
can be identified with the universal cover $\Omega^2 (S^3 \three)$ of
$\Omega^2 S^3$. The resulting map of $\E_2$-spaces 
\[
\Omega^2 (S^3\three)\to \Omega^2S^3
\]
Thomifies to the map $H\Z_{\hat{2}}\to H\F_2$ after composing with the braid representation described in Section \ref{braidrep}. Incidentally, the fibration $\Omega^2 (S^3\three)\to \Omega^2S^3\to S^1$ also leads, as a result of Corollary 4 of \cite{beardsrelative}, to a description of the map $H\Z_{\hat{2}}\to H\F_2$ as a Galois extension with Hopf-Galois object $\sus S^1_+$. 

\subsection{} 

Juxtaposing two braids, each with writhe zero, yields another such braid, defining a (strictly associative, but not braided) monoidal {\bf dewrithed} category
\[
\B^0_\bullet := \coprod_{n \geq 0} [\{n\}/\B^0_n] \;,
\]
with 
\[
\B^0_n := \ker [w_n : \B_n \to \Z] \cong [\B_n,\B_n] \;,
\]
where this last isomorphism arises from noticing that $w_n$, the writhe map, is in fact the abelianization. Thus $\B^m = \{1\}$ if $m = 0,1,2$, and in low degrees we have actions
\[
\B^0_m \times \B^0_n \ni 1 \times \sigma_i \mapsto \sigma_{m+i} \in 
\B^0_{m+n} \;,
\]
\[
\B^0_n \times \B^0_m \ni \sigma_i \times 1 \mapsto \sigma_i \in \B^0_{n+m} \;.
\]
Visualizing braids as displayed vertically and braid multiplication $\beta\cdot\beta'$ as corresponding to placing $\beta'$ below $\beta$, $\B^0_1$ acts on $\B^0_\bullet$ on the left by adding a trivial strand ($\beta \mapsto \beta_{-}$) on the left,  and on the right ($\beta \mapsto \beta_{+}$) by adding a trivial strand on the right. It is not hard to see that if $\beta \in \B^0_n$ then we have compositions 
\[
c_{n,1} \cdot \beta_{+} =  \beta_{-} \cdot c_{n,1} \in \B_{n+1} \;.
\]
For example, if $\beta = \sigma_1 \in \B_2$ then $\beta_{+} = \sigma_1, \; \beta_{-} = \sigma_2 \in \B_3$, while $c_{2,1} = \sigma_1 \cdot \sigma_2 \in \B_3$, yielding 
\[
c_{2,1} \cdot \beta_{+} = \sigma_1 \sigma_2 \sigma_1 = \sigma_2 \sigma_1 \sigma_2 = \beta_{-} \cdot c_{2,1} \;,
\]
The left and right actions thus differ by conjugation with $c_{n,1} \in \B_{n+1}$, i.e.~by the $n$th power of a generator of the group $\Z = \B_{n+1}/\B^0_{n+1}$ of outer automorphisms of $\B^0_{n+1}$.

By construction we have a homotopy pullback diagram:
\[
\xymatrix{
|\B^0_\bullet|\ar@{..>}[r]\ar@{..>}[d] & \mathbb{N}\ar[d]\\
|\B_{\bullet}|\ar[r] & |\B_{\bullet}^\ab|.
}
\] 

Using strengthenings of \cite{mcduffsegal} by \cite{millerpalmer}, we can apply Randal-Williams' methods \cite[Theorem 1.1]{randalwilliamslocal} even though $\B_{\bullet}^0$ is not homotopy commutative: by the remark above, it is enough to use right fractions in Randal William's Theorem 1.1 \cite{randalwilliamslocal}, obtaining a homotopy pullback diagram of group completions:

\[
\xymatrix{
	\Omega B|\B^0_\bullet|\ar@{..>}[r]\ar@{..>}[d] & \mathbb{Z}\ar[d]\\
	\Omega B|\B_{\bullet}|\ar[r] & \Omega B|\B_{\bullet}^\ab|.
}
\] 

The bottom horizontal map in the above diagram is a map of 2-fold loop spaces, but the right vertical map is not, so we cannot say that $\Omega B|\B_{\bullet}^0|$ is an $\E_2$-algebra (and it should not be, because $\B_{\bullet}^0$ is only monoidal). However, by restricting to identity components in the diagram, and using the fact that the inclusion $\ast\hookrightarrow B\Z$ is trivially an infinite loop map (and thus a map of 2-fold loop spaces), we obtain the following homotopy pullback diagram of 2-fold loop spaces:

\[
\xymatrix{
	\Omega B|\B^0_\bullet|_0\ar@{..>}[r]\ar@{..>}[d] & \ast\ar[d]\\
	\Omega^2S^3\ar[r] & B\Z.
}
\] 

We have already shown the fiber of the map $\Omega^2S^3\to B\Z$ to be equivalent to $\Omega^2(S^3\three)$, but now we see that this fiber is further equivalent to the identity component of the group completion of $|\B_{\bullet}^0|$.

\subsection{}  The above analysis admits a geometric description in line with thinking of $\B_{\bullet}$ as being constructed from configuration spaces. Recall from Section \ref{abelianizationmap} that a finite subset $F\subset\C$ induces an ``electrical field,'' determined by a polynomial $p_F(z)=\prod_{\rho\in F}(z-\rho)$ with $F$ as its set of zeros.
The classical discriminant 
\[
\C \supset F \mapsto \Delta(F) := (-1)^{n(n-1)/2}\prod_{\rho \neq \rho' \in 
F}(\rho - \rho') \in \C^\times
\]
($n = \# F$) can be seen to induce the abelianization map $w_n$ on fundamental groups. This gives a geometric model for our map $|\B_{\bullet}|\simeq\Omega^2S^2\to S^1$, with fiber product 
\[
\xymatrix{
	{\bf B}\B^0_n \ar@{.>}[d] \ar@{.>}[r] & \C \ar[d]^\be \\
	B\B_n \ar[r]^{\Delta} & \C^\times }
\]
(where $\be(w) = \exp(2\pi iw)$) homotopy-equivalent to the components of $|\B^0_\bullet|$ (though not equal, so we write $\mathbf{B}$ instead of $B$). In other words, we may think of these components as spaces of configurations $(F,\delta_F)$ augmented by an angular anomaly $\delta_F \equiv (2\pi i)^{-1} \log \Delta(F)$ (mod $\Z$). 

If $F = \{\rho_\alpha \in \C\}$ is an unordered $n$-tuple of distinct points, and $u \in \C^\times$, then $uF := \{ u\rho_\alpha \}$ is another such tuple, with $\Delta(uF) = u^{n(n-1)}\Delta(F)$; consequently $\C$ acts on ${\bf B}\B^0_n$ by 
\[
\Delta(\be(w)F) = \be(\delta_F + n(n-1)w) \;.
\]
Randal-Williams' construction thus provides us with an $\E_1$ model for the group completion of $|\B^0_\bullet|$ as $\Z\times\Omega^2 (S^3 \three)$, with the covering group $\Z \subset \C$ acting by 
\[
k \cdot (F,\delta_F) = (F,\delta_F + n(n-1)k) \;.
\]

\end{document}